\input amstex
\documentstyle{amsppt}
\magnification=\magstep1                        %<====
\hsize6.5truein\vsize8.9truein                  %<====
\NoRunningHeads
\loadeusm

\magnification=\magstep1                        %<====
\hsize6.5truein\vsize8.9truein                  %<====
\NoRunningHeads
\loadeusm

\document

\topmatter

\title
The $L_q$ norm of the Rudin-Shapiro polynomials on subarcs of the unit circle 
\endtitle

\rightheadtext{notes on the Rudin-Shapiro polynomials}

\author
Tam\'as Erd\'elyi
\endauthor

\address Department of Mathematics, Texas A\&M University
College Station, Texas 77843 \endaddress

\email terdelyi\@tamu.edu
\endemail

%\subjclass 11J54,11B83 \endsubjclass

\endthanks

\date June 12, 2023 \enddate


\abstract
Littlewood polynomials are polynomials with each of their coefficients in
$\{-1,1\}$. A sequence of Littlewood polynomials that satisfies a remarkable 
flatness property on the unit circle of the complex plane is given by the
Rudin-Shapiro polynomials.
Let $P_k$ and $Q_k$ denote the Rudin-Shapiro polynomials of degree $n-1$ with $n:=2^k$. 
For polynomials $S$ we define  
$$M_q(S,[\alpha,\beta]) := \left( \frac{1}{\beta-\alpha} 
\int_{\alpha}^{\beta} {\left| S(e^{it}) \right|^q\,dt} \right)^{1/q}\,, \qquad q > 0\,.$$
Let $\gamma := \sin^2(\pi/8)$. We prove that
$$\frac{\gamma}{4\pi}(\gamma n)^{q/2} \leq M_q(P_k,[\alpha,\beta])^q \leq (2n)^{q/2}$$  
for every $q > 0$ and $32\pi/n \leq \beta-\alpha$. The same estimates hold for $P_k$ replaced by $Q_k$.
\endabstract

\endtopmatter

\head 1. Introduction and Notation \endhead
Let $\alpha < \beta$ be real numbers. The Mahler measure $M_{0}(S,[\alpha,\beta])$ is
defined for polynomials $S$ as
$$M_{0}(S,[\alpha,\beta]) := \exp\left(\frac{1}{\beta - \alpha} \int_{\alpha}^{\beta} {\log|S(e^{it})|\,dt} \right)\,.$$
It is well known, see [17] for instance, that
$$M_0(S,[\alpha,\beta]) = \lim_{q \rightarrow 0+}{M_{q}(S,[\alpha,\beta])}\,,$$
where
$$M_{q}(S,[\alpha,\beta]) := 
\left( \frac{1}{\beta-\alpha} \int_{\alpha}^{\beta} {\left| S(e^{it}) \right|^q\,dt} \right)^{1/q}\,, \qquad q > 0\,.$$
It is a simple consequence of the Jensen formula that
$$M_0(S,[0,2\pi]) = |c| \prod_{k=1}^n{\max\{1,|z_k|\}}$$
for every polynomial of the form
$$S(z) = c\prod_{k=1}^n{(z-z_k)}\,, \qquad c,z_k \in {\Bbb C}\,.$$
See [3 p. 271] or [2 p. 3] for instance. 
Let $D := \{z \in {\Bbb C}: |z| < 1\}$ denote the open unit disk of the complex plane.
Let $\partial D :=  \{z \in {\Bbb C}: |z| = 1\}$ denote the unit circle of the complex plane.
Littlewood polynomials are polynomials with each of their coefficients in $\{-1,1\}$. A special
sequence of Littlewood polynomials is the sequence the Rudin-Shapiro polynomials, They appear in Harold
Shapiro's 1951 thesis [21] at MIT and are sometimes called just the Shapiro polynomials.
They also arise independently in Golay's paper [16]. They are remarkably simple to construct
recursively as follows. Let  
$$P_0(z) :=1\,, \qquad Q_0(z) := 1\,,$$
and
$$\split P_{k+1}(z) & := P_k(z) + z^{2^k}Q_k(z)\,, \cr
Q_{k+1}(z) & := P_k(z) - z^{2^k}Q_k(z)\,, \cr \endsplit$$
for $k=0,1,2,\ldots\,.$ Note that both $P_k$ and $Q_k$ are polynomials of degree $n-1$ with $n := 2^k$ 
having each of their coefficients in $\{-1,1\}$. In what follows $P_k$ and $Q_k$ denote the Rudin-Shapiro polynomials 
of degree $n-1$ with $n := 2^k$. 
It is well known, and easy to check by using the parallelogram law, that
$$|P_{k+1}(z)|^2 + |Q_{k+1}(z)|^2 = 2(|P_k(z)|^2 + |Q_k(z)|^2)\,, \qquad z \in \partial D\,.$$
Hence
$$|P_k(z)|^2 + |Q_k(z)|^2 = 2^{k+1} = 2n\,, \qquad z \in \partial D\,. \tag 1.1$$
It is also well known, see Section 4 of [2] or [6] for instance, that 
$$Q_k(z) = (-1)^{k+1}P_k^*(-z)\,, \qquad z \in \partial D\,, \tag 1.2$$
where $P_k^*(z) := z^{n-1}P_k(1/z)$. Hence 
$$|Q_k(z)| = |P_k(-z)|\,, \qquad z \in \partial D\,. \tag 1.3$$
Peter Borwein's book [2] presents a few more basic results on the Rudin-Shapiro 
polynomials. Cyclotomic properties of the Rudin-Shapiro polynomials are discussed in 
[6]. Obviously $M_2(P_k,[0,2\pi]) = 2^{k/2}$ by the Parseval formula. In 1968 Littlewood 
[19] showed that $M_4(P_k,[0,2\pi]) \sim (4^{k+1}/3)^{1/4}$. 
Here, and in what follows, $a_k \sim b_k$ means that $\displaystyle{\lim_{k \rightarrow \infty}{\frac{a_k}{b_k} = 1}}$.
Rudin-Shapiro like polynomials in $L_4$ on the unit circle ${\partial D}$ of the complex plane are studied in [4]. 
Let $K := {\Bbb R} \enskip(\text {mod}\,\, 2\pi)$.
Let $m(A)$ denote the one-dimensional Lebesgue measure of $A \subset K$.
In 1980 Saffari conjectured the following result. He did not publish this conjecture himself,  
and it first appeared in print in the work of Doche and Habsieger [9]. 

\proclaim{Theorem 1.1}
We have 
$$M_q(P_k,[0,2\pi]) = M_q(Q_k,[0,2\pi]) \sim \frac{2^{(k+1)/2}}{(q/2+1)^{1/q}} = \frac{(2n)^{1/2}}{(q/2+1)^{1/q}}$$
for all real exponents $q > 0$. Equivalently, we have 
$$\split & \lim_{k \rightarrow \infty} 
m{\left(\left\{t \in K: \left| \frac{P_k(e^{it})}{\sqrt{2^{k+1}}} \right|^2 \in [\alpha,\beta] \right\}\right)} \cr 
= \, & \lim_{k \rightarrow \infty}
m{\left(\left\{t \in K: \left| \frac{Q_k(e^{it})}{\sqrt{2^{k+1}}} \right|^2 \in [\alpha,\beta] \right\}\right)} 
= 2\pi(\beta - \alpha) \cr \endsplit$$
whenever $0 \leq \alpha < \beta \leq 1$. 
\endproclaim

Theorem 1.1 was proved for all even values of $q \leq 52$ by Doche [8]
and Doche and Habsieger [9]. Rodgers [20] proved Theorem 1.1 for all 
$q > 0$. See also [10]. An application of Theorem 1.1 may be found in [15]. 
An extension of Saffari's conjecture is Montgomery's conjecture below proved 
by Rodgers [20] as well.  

\proclaim{Theorem 1.2}
We have
$$\split & \lim_{k \rightarrow \infty} 
m{\left(\left\{t \in K: \frac{P_k(e^{it})}{\sqrt{2^{k+1}}} \in E \right\}\right)} \cr
= \, & \lim_{k \rightarrow \infty}
m{\left(\left\{t \in K: \frac{Q_k(e^{it})}{\sqrt{2^{k+1}}} \in E \right\}\right)} 
= 2m(E) \cr \endsplit$$
for any rectangle $E \subset D := \{z \in {\Bbb C}: |z| < 1\}\,.$
\endproclaim

In [11] we proved the following lower bound for the Mahler measure of the 
Rudin-Shapiro polynomials on subarcs of the unit circle ${\partial D}$.

\proclaim{Theorem 1.3} 
There is an absolute constant $c > 0$ such that
$$M_0(P_k,[\alpha,\beta]) \geq cn^{1/2}$$ 
for all $k \in {\Bbb N}$ and for all $\alpha, \beta \in {\Bbb R}$ such that
$$\frac{32\pi}{n} \leq \frac{(\log n)^{3/2}}{n^{1/2}} \leq \beta - \alpha \leq 2\pi\,.$$
The same lower bound holds for $M_0(P_k,[\alpha,\beta])$ replaced by $M_0(Q_k,[\alpha,\beta])$.
\endproclaim

It looks plausible that Theorem 1.3 holds whenever $32\pi/n \leq \beta - \alpha\,,$
but we have not been able to handle the case $32\pi/n \leq \beta - \alpha \leq (\log n)^{3/2}n^{-1/2}$. 
Nevertheless our Theorem 2.2 gives a lower bound for the values 
$M_q(P_k,[\alpha,\beta])$ and $M_q(Q_k,[\alpha,\beta])$ for every $q > 0$ and $32\pi/n \leq \beta - \alpha$.  
See also [7] on sums of monomials with large Mahler measure on subarcs of the unit circle ${\partial D}$.
In [13] the asymptotic values of $M_0(P_k,[0,2\pi])$ and $M_0(P_k,[0,2\pi])$, conjectured by Saffari, 
have been found. Namely in [13] we showed the following. 

\proclaim{Theorem 1.4}
We have 
$$\lim_{n \rightarrow \infty}{\frac{M_0(P_k,[0,2\pi])} {n^{1/2}}} 
= \lim_{n \rightarrow \infty}{\frac{M_0(Q_k,[0,2\pi])}{n^{1/2}}}
= \left( \frac 2e \right)^{1/2}\,.$$  
\endproclaim

Properties of the Rudin Shapiro polynomials have played a a central role in [1] 
as well as in [14] to prove a longstanding conjecture of Littlewood on the existence of flat 
Littlewood polynomials $S_n$ of degree $n$ satisfying the inequalities
$$c_1 n^{1/2} \leq |S_n(e^{it})| \leq c_2 n^{1/2}, \qquad t \in {\Bbb R} \,,$$
with absolute constants $c_1 > 0$ and $c_2 > 0$.

\head New Results \endhead
Let $\gamma := \sin^2(\pi/8)$ and $n := 2^k$ The Lebesgue measure of a set $E \subset {\Bbb R}$ is denoted by $m(E)$.

\proclaim{Theorem 2.1}
Let $E :=  \{t \in [\alpha,\beta]: |P_k(t)| \geq \gamma n\}\,.$
We have 
$$m(E) \geq \frac{(\beta - \alpha)\gamma}{4\pi}$$  
for every $32\pi/n \leq \beta - \alpha$. The same estimate holds for $P_k$ replaced by $Q_k$.
\endproclaim

\proclaim{Theorem 2.2}
We have 
$$\frac{\gamma}{4\pi} (\gamma n)^{q/2} \leq M_q(P_k,[\alpha,\beta])^q  \leq (2n)^{q/2}$$
for every $q > 0$ and $32\pi/n \leq \beta - \alpha$. The same estimate holds for $Q_k$.
\endproclaim

\head 3. Lemmas \endhead 
Let $n:=2^k$, $\gamma := \sin^2(\pi/8)$, $z_j := e^{it_j}$, $t_j := 2\pi j/n$, $j \in {\Bbb Z}$.
 
\proclaim{Lemma 3.1}
We have
$$\max \{|P_k(z_j)|^2,|P_k(z_{j+r})|^2\} \geq \gamma 2^{k+1} = 2\gamma n\,, \quad r \in \{-1,1\}\,,$$
for every $j=2u$, $u \in {\Bbb Z}$. The same estimate holds for $P_k$ replaced by $Q_k$.
\endproclaim

Lemma 3.1 tells us that the modulus of the Rudin-Shapiro polynomials $P_k$ is certainly not smaller than $(2\gamma n)^{1/2}$ 
at least at one of any two consecutive $n$-th root of unity, where $n := 2^k$. This is a crucial observation proved in [11] 
and plays a key role in [12], [13], [14] and [15] as well. Our Lemma 3.2 below follows from Lemma 3.1 
reasonably simply.

\proclaim{Lemma 3.2}
We have
$$|P_k(e^{it})|^2 \geq \gamma n\,, \qquad t \in [t_j - \gamma/n, t_j + \gamma/n]\,,$$
for every $j \in {\Bbb Z}$ such that
$$|P_k(z_j)|^2 \geq \gamma 2^{k+1} = 2\gamma n\,. \tag 3.1$$  
The same estimate holds with $P_k$ replaced by $Q_k$.
\endproclaim

\demo{Proof of Lemma 3.2}
By (1.3) it is sufficient to prove the lemma only for $P_k$.
The proof of the lemma is a simple combination of the Mean Value Theorem and Bernstein's inequality 
applied to the nonnegative trigonometric polynomial $R_k$ of degree $n-1$ with $n=2^k$ defined by 
$R_k(t) := P_k(e^{it})P_k(e^{-it})$. Recall that (1.1) implies $0 \leq R_k(t) = |P_k(e^{it})|^2 \leq 2n$ 
for every $t \in {\Bbb R}$. Note also that the Bernstein factor is $n/2$ rather than $n$ for the class of 
nonnegative trigonometric polynomials of degree at most $n$, see Lemma 3.3 below. Suppose 
$j \in {\Bbb Z}$ satisfies (3.1) and $t \in {\Bbb R}$ satisfies $|t-t_j| \leq \gamma/n$. Then by the Mean Value Theorem 
there is a $\xi$ between $t_j$ and $t$ such that 
$$R_k(t_j) - R_k(t) \leq |R_k(t_j) - R_k(t)| = |t_j-t||R_k^\prime(\xi)| 
\leq \frac{\gamma}{n} \frac n2 \max_{\tau \in K}{\{R_k(\tau)} \leq \frac{\gamma}{n} \frac n2 2n = \gamma n\,.$$
Therefore, recalling (3.1), we get
$$R_k(t) \geq R_k(t_j) - \gamma n = 2\gamma n - \gamma n = \gamma n\,, \qquad t \in [t_j-\gamma/n\,, t_j+\gamma/n]\,.$$
\qed \enddemo

Let $K := {\Bbb R} \enskip(\text {mod}\,\, 2\pi)$, as before. 

\proclaim{Lemma 3.3}
We have 
$$\max_{\tau \in K}{|T^\prime(\tau)|} \leq \frac n2 \max_{\tau \in K}{T(\tau)}$$
for every trigonometric polynomial $T$ of degree at most $n$ that is nonnegative on ${\Bbb R}$. 
\endproclaim

\demo{Proof of Lemma 3.3}
Suppose $T$ is a trigonometric polynomial of degree at most $n$ that is nonnegative on ${\Bbb R}$. 
The Bernstein inequality, see [3] for instance, asserts that 
$$\max_{\tau \in K}{||Q^\prime(\tau)|} \leq n \max_{\tau \in K}{|Q(\tau)|}$$ 
for every real trigonometric polynomial $Q$ of degree at most $n$.
Applying the Bernstein inequality to the real trigonometric polynomial 
$Q := T - M$ of degree at most $n$ with $\displaystyle{M := \frac 12 \max_{\tau \in K}{|Q(\tau)|}}$
gives the lemma.
\qed \enddemo

\head 4. Proof of the theorems \endhead

\demo{Proof of Theorem 2.1}
By (1.3) it is sufficient to prove the theorem only for $P_k$. 
Observe that Lemmas 3.1 and 3.2 imply that $E$ contains at least $\displaystyle{\frac{\beta - \alpha)n}{4\pi}-4}$ disjoint intervals 
of length at least $2\gamma/n$, hence 
$$m(E) \geq \left( \frac{(\beta - \alpha)n}{4\pi}-4 \right) \frac{2\gamma}{n} \geq \frac{\beta-\alpha}{8\pi} \frac{2\gamma}{n} = \frac{(\beta - \alpha)\gamma}{4\pi}$$  
whenever $32\pi/n \leq \beta - \alpha$.
\qed \enddemo

\demo{Proof of Theorem 2.2}  
By (1.2) it is sufficient to prove Theorem 2.1 for $P_k$. The upper bound of the theorem follows immediately from (1.1). 
Now we prove the lower bound of the theorem. Using Theorem 2.1 we have 
$$\split M_q(P_k,[\alpha,\beta])^q & := 
\frac{1}{\beta-\alpha} \int_{\alpha}^{\beta} {|P_k(t)|^q\,dt} \geq \frac{1}{\beta-\alpha}\int_E{|P_k(t)|^q \, dt} \cr 
& \geq \frac{1}{\beta-\alpha} m(E)(\gamma n)^{q/2} \geq \frac{1}{\beta-\alpha}\frac{(\beta - \alpha)\gamma}{4\pi} (\gamma n)^{q/2}\cr 
& \geq \frac{\gamma}{4\pi}(\gamma n)^{q/2} \cr \endsplit$$ 
whenever $32\pi/n \leq \beta-\alpha$.
\qed \enddemo

\head 5. More observations and problems \endhead

Let $P_k$ and $Q_k$ be the usual Rudin-Shapiro polynomials of degree $n-1$ with $n:=2^k$.

As for $k \geq 1$ both $P_k$ and $Q_k$ have odd degree $n-1 = 2^k-1$, both $P_k$ and $Q_k$ have at least one
real zero. The fact that for $k \geq 1$ both $P_k$ and $Q_k$ have exactly one real zero was proved
by Brillhart in [5]. Another interesting observation made in [6] is the fact that $P_k$ and $Q_k$ 
cannot vanish at any roots of unity different from $-1$ and $1$. In [12] we proved that the Rudin-Shapiro 
polynomials $P_k$ and $Q_k$  have only $o(n)$ zeros on the unit circle ${\partial D}$.  Observe, 
see [6] for instance, that 
$$P_k(1) = 2^{[(k+1)/2]}, \qquad Q_k(-1) = (-1)^{k+1} 2^{[(k+1)/2]},$$
and $$P_k(-1) = Q_k(1) = \frac 12(1 + (-1)^k)2^{[k/2]},$$
where $[x]$ denotes the integer part of a real number $x$.

\proclaim{Problem 5.1}
Is it true that if $k$ is odd then $P_k$ has a zero on the unit circle ${partial D}$ only at $-1$ and
$Q_k$ has a zero on the unit circle ${\partial D}$ only at $1$, while if $k$ is even then neither $P_k$ nor $Q_k$ 
has a zero on the unit circle ${\partial D}$?
\endproclaim
Combining (1.2) with the observation that the Rudin-Shapiro polynomials $P_k$ and $Q_k$ of degree $n-1$ with 
$n := 2^k$ have only $o(n)$ zeros on the unit circle ${\partial D}$, we can deduce that the products 
$P_kQ_k$ have $n-o(n)$ zeros in the open unit disk $D$, where $o(n)$ denotes real numbers 
such that $o(n)/n$ converges to $0$ as $n$ tends to $\infty$.   

\proclaim{Problem 5.2}
Is there an absolute constant $c > 0$ such that both of the Rudin-Shapiro polynomials $P_k$ and $Q_k$  
have at least $cn$  zeros in the open unit disk $D$?    
\endproclaim

\proclaim{Problem 5.3}
Is it true that both of the Rudin-Shapiro polynomials $P_k$ and $Q_k$ have $n/2-o(n)$, zeros in the open unit disk $D$?
\endproclaim

\proclaim{Problem 5.4}
Is it true that Theorem 1.3 remains valid for all $32\pi/n \leq \beta-\alpha \leq  2\pi$? 
\endproclaim

\proclaim{Problem 5.5}
Is there an absolute constant $c > 0$ such that 
$$M_0(|P_k|^2-n,[0,2\pi]) := 
\exp\left(\frac{1}{2\pi} \int_0^{2\pi} {\log||P_k(e^{it})|^2-n|\,dt} \right)\geq c n^{1/2}?$$ 

\endproclaim

\head 6. A connection to sew-reciprocal polynomials \endhead

A polynomial $S$ of the form
$$S(z) = \sum_{j=0}^{2m}{a_jz^j}\,, \qquad a_j \in {\Bbb R}\,, \quad a_{2m} \neq 0\,,$$
is called skew-reciprocal if
$$a_{m-j} = (-1)^ja_{m+j}\,, \qquad j=1,2,\ldots,m\,. \tag 6.2)$$ 

A beautiful observation of Mercer [18] states the following.

\proclaim{Theorem 6.1} 
Skew-reciprocal Littlewood polynomials do not have any zeros on the unit circle ${\partial D}$.
\endproclaim 

The Rudin-Shapiro polynomials $P_k$ and $Q_k$ of degree $n-1$ with $n := 2^k$  
are quite close to be skew-reciprocal. However, as the degrees of $P_k$ and $Q_k$ are odd, 
Theorem 6.1 does not apply to the Rudin-Shapiro polynomials. Having a middle term in the 
polynomial $S$ in the proof below is crucial.   

\demo{Proof of Theorem 6.1}
Let $S$ be a skew-reciprocal Littlewood polynomial of the form 
$$S(z) = \sum_{j=0}^{2m}{a_jz^j}, \quad a_j \in \{-1,1\}\,, \quad j = 0,1, \ldots, 2m\,, \quad a_{2m} \neq 0\,,$$ 
with
$$a_{m-j} = (-1)^j a_{m+j}\,, \qquad j=1,2,\ldots,m\,.$$ 
For notational convenience we assume that $m=2\mu$ is even; the proof in the case when $m=2\mu-1$ is odd can be 
handled similarly. 
We have $z^{-m}S(z) = A(z) + B(z)$, where the function
$$A(z) := \sum_{j=0}^{\mu}{a_{m+2j}}{(z^{2j} + z^{-2j})}\,, \qquad z \in \partial D\,,$$ 
takes purely real values on the unit circle ${\partial D}$, and the function 
$$B(z) := \sum_{j=1}^{\mu}{a_{m+2j-1}(z^{2j-1} - z^{-2j-1})}\,, \qquad z \in \partial D\,,$$ 
takes purely imaginary values on the unit circle ${\partial D}$. Suppose to the contrary that $S$ vanishes at 
a point $z_0$ on the unit circle $\partial D$. Then $z_0$ is a common zero of $A$ and $B$. 
We study the greatest common divisor of the polynomials $\widetilde{A}(z) := z^m A(z)$ and $\widetilde{B}(z) := z^mB(z)$ 
over the field ${\bold F}_2$. We have 
$$\widetilde{A}(z) -z\widetilde{B}(z) = \sum_{j=0}^m{z^{2j}} - z\sum_{j=1}^m{z^{2j-1}} = 1$$
over the field ${\bold F}_2$, showing that the greatest common divisor of the polynomials $\widetilde{A}$ and $\widetilde{B}$ 
over the field ${\bold F}_2$ is $1$. Hence $A(z)$ and $B(z)$ cannot have a common zero on the unit circle ${\partial D}$, 
a contradiction.
\qed \enddemo 

Note that the same approach works to prove that skew-reciprocal polynomials with only odd coefficients do not have any 
zeros on the unit circle ${\partial D}$.

\enddocument